\documentclass{amsart}
\usepackage{amsfonts,amssymb, amsmath,amsthm,amscd} 
\usepackage{fullpage}
\usepackage{url}
\usepackage{enumerate}
\usepackage{lmodern}
\newtheorem{theorem}{Theorem}[section]

\newtheorem{proposition}[theorem]{Proposition}
\newtheorem{corollary}[theorem]{Corollary}

\theoremstyle{definition}

\newenvironment{definition}[1][Definition]{\begin{trivlist}
\item[\hskip \labelsep {\bfseries #1}]}{\end{trivlist}}
\newenvironment{remark}[1][Remark]{\begin{trivlist}
\item[\hskip \labelsep {\bfseries #1}]}{\end{trivlist}}

\DeclareFontEncoding{OT2}{}{}
\newcommand{\textcyr}[1]{{\fontencoding{OT2}\fontfamily{wncyr}\fontseries{m}\fontshape{n}\selectfont #1}}

\author{Jeanine Van Order}
\address{Fakult\"at f\"ur Mathematik, Universit\"at Bielefeld}
\email{jvanorder@math.uni-bielefeld.de}
\subjclass{Primary 11F67; Secondary 11F25, 11G40, 11R23}
\begin{document}

\title{Rankin-Selberg $L$-functions in cyclotomic towers, II}

\begin{abstract} 

Following the prequel work \cite{VO3}, we prove a generalization of ``Mazur's conjecture" for $L$-functions 
of elliptic curves in abelian extensions of imaginary quadratic fields, including the assertion that the Mordell-Weil 
rank of an elliptic curve in the ${\bf{Z}}_p^2$-extension is finitely generated modulo Heegner points. 
The novelty of the approach here is to use the existence of a suitable $p$-adic $L$-function to reduce the problem 
to a minimal nonvanishing criterion, which should be applicable to a broader class of problems than considered here. \\


\end{abstract}

\maketitle

\section{Introduction}

Let $E$ be an elliptic curve of conductor $N$ defined over the rational number field. We know that $E$ is modular by fundamental work of Wiles \cite{Wi}, 
Taylor-Wiles \cite{TW} and Breuil-Conrad-Diamond-Tayor \cite{BCDT}, hence parametrized by a cuspidal newform $f$ of weight $2$, level $N$ and trivial 
nebentype character. Let $K$ be an imaginary quadratic extension of discriminant $D < 0$ prime to $N$. Fix an odd prime $p$ which does not divide $N$. 
Let $K_{\infty}$ denote the ${\bf{Z}}_p^2$-extension of $K$, with Galois group $G = \operatorname{Gal}(K_{\infty}/K) \approx {\bf{Z}}_p^2$. 
Hence, $K_{\infty}$ is a pro-$p$-extension of $K$ contained in the maximal abelian extension $R_{\infty}$ of $K$ unramified outside of $p$.
This latter abelian extension $R_{\infty}$ can be described more explicitly as the compositum $K[p^{\infty}]K(\mu_{p^{\infty}})$, 
where $K[p^{\infty}] = \bigcup_{n \geq 0}K[p^n]$ is the union of all ring class extensions of $p$-power conductor of $K$, 
and $K(\zeta_{p^{\infty}}) = \bigcup_{n \geq 0} K(\zeta_{p^n})$ the extension obtained by adjoining all $p$-power roots of unity to $K$. 
Let $\mathcal{G} = \operatorname{Gal}(R_{\infty}/K)$, $\varOmega = \operatorname{Gal}(K[p^{\infty}]/K)$ and 
$\varGamma = \operatorname{Gal}(K(\zeta_{p^{\infty}})/K)$ denote the corresponding profinite Galois groups.
The torsion subgroup $G_0 = \mathcal{G}_{\operatorname{tors}}$ is finite, 
and the quotient $G = \mathcal{G} \operatorname{mod} G_0 \approx {\bf{Z}}_p^2$ determines the extension $K_{\infty}/K$.
More precisely, the complex conjugation automorphism in $\operatorname{Gal}(K/{\bf{Q}})$ acts by conjugation on $G$, and the respective $\pm 1$-
eigenspaces define
${\bf{Z}}_p$-quotients $\Omega = \operatorname{Gal}(D_{\infty}/K)$ 
of $\varOmega$ and $\Gamma = \operatorname{Gal}(K^{\operatorname{cyc}}/K)$ of $
\varGamma$ with 
$G \approx \Omega \times \Gamma \approx {\bf{Z}}_p^2$. 

Let ${\bf{Z}}_p[[\mathcal{G}]]$ denote the ${\bf{Z}}_p$-Iwasawa algebra of a profinite group $\mathcal{G}$. Consider the $p$-primary Selmer group 
$\operatorname{Sel}(E/K_{\infty})$ of $E$ over $K_{\infty}/K$, which fits into the exact descent sequence of discrete ${\bf{Z}}_p[[G]]$-modules 
\begin{align}\label{SES} 0 \longrightarrow E(K_{\infty}) \otimes {\bf{Q}}_p/{\bf{Z}}_p \longrightarrow \operatorname{Sel}(E/K_{\infty}) \longrightarrow 
\mbox{\textcyr{Sh}}(E/K_{\infty})(p) \longrightarrow 0. \end{align} 
Here, $E(K_{\infty})$ denotes the Mordell-Weil group of $E$ over $K_{\infty}$, and $\mbox{\textcyr{Sh}}(E/K_{\infty})(p)$ the $p$-primary subgroup of the 
Tate-Shafarevich group of $E$ over $K_{\infty}$. The Pontryagin dual 
$X(E/K_{\infty}) = \operatorname{Hom} (\operatorname{Sel}(E/K_{\infty}), {\bf{Q}}_p/{\bf{Z}}_p)$ 
has the structure of a compact, finitely-generated ${\bf{Z}}_p[[G]]$-module. A deep theorem of Kato \cite{Ka}, using Rohrlich's nonvanishing 
theorem \cite{Ro}, shows that the dual Selmer group $X(E/ K^{\operatorname{cyc}})$ has the structure of a finitely-generated torsion 
${\bf{Z}}_p[[\Gamma]]$-module. Hence, $X(E/ K^{\operatorname{cyc}})$ has a ${\bf{Z}}_p[[\Gamma]]$-characteristic power series 
$g(E/K^{\operatorname{cyc}}) \in {\bf{Z}}_p[[\Gamma]]$, which after fixing a topological generator $\gamma \in \Gamma$, 
we view as a power series in the indeterminate $T = \gamma - 1$ 
as $g(E/K^{\operatorname{cyc}})(T) \in {\bf{Z}}_p[[T]] = {\bf{Z}}_p[[\gamma -1]] \approx {\bf{Z}}_p[[\Gamma]]$. 
Starting with this result, we show in \cite{VO3} that the two-variable dual Selmer group $X(E/K_{\infty})$ has the  
structure of a finitely generated torsion ${\bf{Z}}_p[[G]]$-module, and hence has a two-variable characteristic power series 
$g(E/K_{\infty}) \in {\bf{Z}}_p[[G]] \approx {\bf{Z}}_p[[T_1, T_2]]$. Moreover, if the exponent of $p$ in the cyclotomic 
power series $g(E/K^{\operatorname{cyc}})(T)$ vanishes, i.e.~if the cyclotomic $\mu$-invariant of $X(E/K^{\operatorname{cyc}})$ vanishes, 
then we even have a ${\bf{Z}}_p[[\operatorname{Gal}(K_{\infty}/K^{\operatorname{cyc}})]]$-module isomorphism 
\begin{align}\label{IMS} X(E/K_{\infty}) &\longrightarrow {\bf{Z}}_p[[\operatorname{Gal}(K_{\infty}/K^{\operatorname{cyc}})]]^{\lambda_E(K)}, \end{align}
where $\lambda_E(K)$ denotes the Weierstrass degree of $g(E/K^{\operatorname{cyc}})(T)$, i.e.~the cyclotomic $\lambda$-invariant.
It is possible to deduce from the nonvanishing theorems of Cornut-Vatsal \cite[Theorems 1.4-5]{CV} for totally real fields 
that $E(K_{\infty})$ is typically finitely generated (see \cite[Proposition 3.14, (12)]{VO3}, and note that $p \nmid h_K$ is not needed). 
In brief, writing $\epsilon \in \lbrace \pm 1 \rbrace$ to denote the sign of the Hasse-Weil $L$-function $L(E/K, s)$, 
\begin{align}\label{cork} \operatorname{corank}_{{\bf{Z}}_p[[\operatorname{Gal}(K_{\infty}/K^{\operatorname{cyc}})]]} 
\left( E(K_{\infty}) \otimes {\bf{Q}}_p/{\bf{Z}}_p \right) 
&= \begin{cases} 0 &\text{ if $\epsilon = +1$}\\ 1 &\text{ if $\epsilon = -1$}\end{cases}. \end{align} 
Assuming the cyclotomic $\mu$-invariant vanishes, it can then be deduced from $(\ref{IMS})$ that
\begin{align*} \operatorname{corank}_{{\bf{Z}}_p[[\operatorname{Gal}(K_{\infty}/K^{\operatorname{cyc}})]]} \left( \mbox{\textcyr{Sh}}(E/K_{\infty})(p) \right) 
&= \begin{cases} \lambda_E(K) &\text{ if $\epsilon = +1$}\\ \lambda_E(K)-1 &\text{ if $\epsilon = -1$.}\end{cases} \end{align*}

It is also easy to deduce from $(\ref{cork})$ that $E(K_{\infty})$ must be finitely generated if $\epsilon = +1$, 
and that $E(K_{\infty})/E(D_{\infty})$ must be finitely generated if $\epsilon = -1$ (although it is not stated explicitly in the introduction of \cite{VO3}):

\begin{theorem}\label{MWr} 

Let $E$ be an elliptic curve of conductor $N$ defined over the rationals. Let $K$ be an imaginary quadratic field of discriminant $D$ prime to $N$. 
Let $p > 2$ be a rational prime. Assume that $p$ does not divide $D N$, and the eigenform $f$ parametrizing $E$ is $p$-ordinary, 
i.e. that the eigenvalue of the Hecke operator acting on $f$ is a $p$-adic unit. If $\epsilon = +1$, 
then $E(K_{\infty})$ is finitely generated; if $\epsilon = -1$, then $E(K_{\infty})/E(D_{\infty})$ is finitely generated. \end{theorem}

In this work, we consider the corresponding nonvanishing behaviour of the two-variable $p$-adic $L$-function 
$\mathcal{L}_p(E/K_{\infty}) = \mathcal{L}_p(f/K_{\infty}) \in {\bf{Z}}_p[[G]]$ constructed by Hida \cite{Hi} and Perrin-Riou \cite{PR88}. 
We shall show how a similar style of argument as used to derive Theorem \ref{MWr} can be applied to this 
$p$-adic $L$-function $\mathcal{L}_p(f/K_{\infty})$ to derive a strong nonvanishing result for the corresponding family
of complex (degree-four, non self-dual) central values without using any analysis (cf.~e.g.~\cite{Ro} or \cite{VO5}). 
In particular, using a purely algebraic argument, we derive a minimal complex nonvanishing criterion (Proposition \ref{criterion} below) 
required to establish such results. Putting this together with the nonvanishing estimates of the prequel work \cite{VO5} then allows us to derive 
our main theorem below, which can be viewed as an asymptotic generalization of Mazur's conjecture for anticyclotomic extensions. 

In brief, the two-variable $p$-adic $L$-function $\mathcal{L}_p(f/K_{\infty})$ we refer to is a ${\bf{Z}}_p$-valued measure on $G$ which interpolates 
the complex central values $L(E/K, \mathcal{W}, 1) \approx L(1/2, f \times \mathcal{W})$, 
for $\mathcal{W}$ any nontrivial finite-order Hecke character of $K$ corresponding under the 
reciprocity map to finite-order characters of $G$. The so-called two-variable main conjecture of Iwasawa theory posits
that the principal ideal generated by this $\mathcal{L}_p(E/K_{\infty})$ equals that generated by $g(E/K_{\infty})$ in the 
Iwasawa algebra ${\bf{Z}}_p[[G]]$: $\left( \mathcal{L}_p(E/K_{\infty}) \right) = \left( g(E/K_{\infty}) \right)$ in ${\bf{Z}}_p[[G]]$. 
In this direction, we know the following:

\begin{theorem} Let $E/{\bf{Q}}$ be an elliptic curve of conductor $N$ as above, where the factorization of $N$ in $K$ is given as 
$N = N^+ N^{-}$ for $N^{+}$ divisible only by primes that split in $K$, and $N^{-}$ divisible only by primes that remain inert in $K$. 
Assume (i) $N^{-}$ is the squarefree product of an odd number of primes (whence $\epsilon = +1$), and (ii)  the residual 
Galois representation $\overline{\rho}_E$ associated to $E$ is surjective. Then, taking $p \nmid ND$ as above to be a prime of 
good ordinary reduction for $E$, the two-variable main conjecture is true:
\begin{align}\label{2VMC} \left( \mathcal{L}_p(E/K_{\infty}) \right) = \left(g(E/K_{\infty}) \right) \text{as principal ideals in } {\bf{Z}}_p[[G]] \end{align} \end{theorem}

\begin{proof} See Skinner-Urban \cite[Theorem 3.30]{SU}.
The result can also be deduced from \cite[Proposition 1.7]{VO3}, using as input the anticyclotomic main conjecture derived in 
Skinner-Urban \cite[Theorem 3.3.7]{SU}, which in turn uses as input the Euler system construction of Bertolini-Darmon \cite{BD}. 
To invoke \cite[Proposition 1.7]{VO3}, we use generalizations of this latter construction to totally real number 
fields which factor through the cyclotomic ${\bf{Z}}_p$-extension of ${\bf{Q}}$, as given to prove \cite[Theorem 1.2]{VO2} (cf.~\cite{Lo} and \cite{Nek}). 
However, the additional hypotheses \cite[Theorem 1.2, (B) (C)]{VO2} then need to be assumed, these being a multiplicity one hypothesis for 
certain character groups associated to the residual Galois representation attached to $E$ and Ihara's lemma for Shimura curves (over totally real fields) respectively. \end{proof}
 
Now, we use the structure of this $p$-adic $L$-function $\mathcal{L}_p(f/K_{\infty})$
to generalize the nonvanishing theorems of Rohrlich \cite{Ro}, \cite{Ro2}, Vatsal \cite{Va}, and Cornut \cite{Cor}
for the central values $L(E/K, \mathcal{W}, 1) \approx L(1/2, f \times \mathcal{W})$ as follows. 
Let us decompose any finite-order character $\mathcal{W}$ of $\mathcal{G} = \operatorname{Gal}(R_{\infty}/K)$ into its dihedral and cyclotomic parts as 
$\mathcal{W} = \rho \psi = \rho \chi \circ {\bf{N}}$, where $\rho$ is a character of the dihedral Galois group $\varOmega = \operatorname{Gal}(K[p^{\infty}]/K)$
and $\psi = \chi \circ {\bf{N}}$ is a character of the cyclotomic Galois group $\varGamma = \operatorname{Gal}(K(\mu_{p^{\infty}})/K)$. 
We shall always identify such characters with their corresponding Hecke characters of $K$, 
so that $\rho$ is a ring class character of $K$ of some $p$-power conductor, and $\psi$ arises from a Dirichlet
character $\chi$ of some $p$-power conductor after composition with the norm homomorphism ${\bf{N}}$. Note that 
a classical construction due to Hecke associates to any such character $\mathcal{W}$ a theta series $\Theta(\mathcal{W})$, 
which is a modular form of weight $1$, level $\Delta  = \vert D \vert c(\mathcal{W})$, and nebentype character 
$\omega \mathcal{W}\vert_{\bf{Q}} = \omega\chi^2$, where $\omega = \omega_D$ is the Dirichlet character associated to the quadratic extension $K$. 
We consider the Rankin-Selberg $L$-function of $f$ times any such theta series $\Theta(\mathcal{W})$, normalized to have central value at $s = 1/2$, 
which we denote here by the symbol $L(s, f \times \mathcal{W})$. It is well-known that this $L$-function has an analytic continuation, 
and that its completed $L$-function \begin{align*} \Lambda(s, f \times \mathcal{W}) &= (N\Delta)^s 
\Gamma_{\bf{R}}(s+1/2)\Gamma_{\bf{R}}(s+3/2)L(s, f \times \mathcal{W}) \end{align*} 
satisfies the functional equation \begin{align}\label{FE} \Lambda(s, f \times \mathcal{W}) &= 
\epsilon(s, f \times \mathcal{W})\Lambda(1-s, f \times \overline{\mathcal{W}}). \end{align} 
Here, $\Gamma_{\bf{R}}(s)$ denotes $\pi^{-s/2}\Gamma(s/2)$, with $\epsilon(s, f \times \mathcal{W})$ the epsilon 
factor of $\Lambda(s, f \times \mathcal{W})$, and $\overline{\mathcal{W}} = \mathcal{W}^{-1}$ the contragredient character associated to $\mathcal{W}$. 
The epsilon factor at the central point $\epsilon(1/2, f \times \mathcal{W})$ 
is a complex number of modulus one known as the {\it{root number}} of $L(s, f \times \mathcal{W})$. 
In this particular setting, assuming that the level $N$ of $f$ is prime to the discriminant 
$D$ of $K$, one can see from the classical derivation of $(\ref{FE})$ via convolution that 
the root number $\epsilon(1/2, f \times \mathcal{W})$ is given by the formula 
\begin{align}\label{RN} \epsilon(1/2, f \times \mathcal{W}) &= \mathfrak{K}_{\omega\chi^2} \cdot \omega\chi^2(- N), \end{align} 
where $\mathfrak{K}_{\omega\chi^2}$ is the complex number of modulus one defined by 
\begin{align}\label{GS} \mathfrak{K}_{\omega\chi^2} &= \tau(\omega\chi^2)^4/ c(\omega\chi^2)^2. 
\end{align} Here, $\tau(\omega\chi^2)$ is the Gauss sum of the Dirichlet character $\omega\chi^2$ 
(see \cite{Li} or \cite{JL}), and $c(\omega \chi^2)$ the conductor of $\omega \chi^2$.
The $L$-function $L(s, f \times \mathcal{W})$ is said to be self-dual if the coefficients in its Dirichlet series expansion are real-valued, 
in which case its root number $\epsilon(1/2, f \times \mathcal{W})$ takes values in the set $\lbrace \pm 1 \rbrace$. 
This is the case when $\mathcal{W} = \rho$ is a ring class character of $K$. In fact, when $\mathcal{W} = \rho$ is such a character, it is easy to see that the functional 
equation $(\ref{FE})$ relates the same completed $L$-function on either side, i.e.~that  $\Lambda(s, f \times \mathcal{W}) = \Lambda(s, f \times \overline{\mathcal{W}})$. 
This is a consequence of the well-known fact that ring class characters are equivariant with respect to complex conjugation after fixing an embedding 
of $K$ into ${\bf{C}}$ (cf.~\cite[p.~384]{Ro2}). Thus in this particular setting with $\mathcal{W} = \rho$, the condition that the root number $\epsilon(1/2, f \times \mathcal{W})$ 
equal $-1$ implies the vanishing of the associated central value $L(1/2, f \times \mathcal{W})$. To distinguish this particular case in all that follows, 
we make the following definition: We define a pair $(f, \mathcal{W})$ to be {\it{exceptional}} if (i) $\mathcal{W} = \rho$ is a ring class character 
and (ii) the root number $\epsilon(1/2, f \times \mathcal{W})$ is $-1$. We then define a pair $(f, \mathcal{W})$ to be {\it{generic}} if it is not exceptional in this sense. 

Implicitly, we shall use the following theorem of Shimura \cite{Sh2} (cf.~\cite{Ro} and \cite{Va}).
Let $F$ denote the Hecke field of $f$, the finite extension of ${\bf{Q}}$ obtained by adjoining 
the Hecke eigenvalues of $f$. Let $F(\mathcal{W})$ denote the cyclotomic extension of $F$ obtained by adjoining 
the values of $\mathcal{W}$. Shimura \cite[Theorem 4]{Sh} shows that the values 
\begin{align}\label{value} \mathcal{L}(1/2, f \times \mathcal{W}) &= L(1/2, f \times \mathcal{W})/\langle f, f \rangle \end{align} 
lie in the number field $F(\mathcal{W})$, where $\langle f, f \rangle$ denotes the Petersson inner product of $f$ with itself. 
In particular, these values $(\ref{value})$ are algebraic. Given a Hecke character $\mathcal{W}$ of $K$ and an 
embedding $F(\mathcal{W}) \rightarrow {\bf{C}}$, let $\mathcal{W}^{\sigma}$ denote the character defined on 
ideals $\mathfrak{a}$ of $\mathcal{O}_K$ by the rule $\mathfrak{a} \mapsto \mathcal{W}(\mathfrak{a})^{\sigma}$. Shimura \cite[Theorem 4]
{Sh} in fact shows that an automorphism $\sigma$ of $F(\mathcal{W})$ acts on the value 
$\mathcal{L}(1/2, f \times \mathcal{W})$ by the rule \begin{align*} \mathcal{L}(1/2, f \times \mathcal{W})^{\sigma} &= \mathcal{L}(1/2, f^{\tau} \times \mathcal{W}^{\sigma}),\end{align*} 
where $\tau$ denotes the restriction of $\sigma$ to $F$. Thus, writing $P_{[\mathcal{W}]}$ to denote the set of characters $\mathcal{W}^{\sigma}$ for $\sigma$ 
ranging over embeddings $F(\mathcal{W}) \rightarrow {\bf{C}}$ fixing $F$, we see that the set of values $\mathcal{L}(1/2, f \times \mathcal{W}^{\sigma})$ are Galois conjugate. 
In particular, it follows that if  $L(1/2, f \times \mathcal{W})$ vanishes for some character $\mathcal{W} \in P_{[\mathcal{W}]}$, then the  
same is true for {\it{all}} characters in the set $P_{[\mathcal{W}]}$. Equipped with this notion of Galois conjugacy,
we show the following result, where the special case of $k=1$ is due to Cornut \cite{Cor}. 

\begin{theorem}\label{main} 

Let $f \in S_2(\Gamma_0(N))$ be a cuspidal eigenform of weight $2$ on $\Gamma_0(N)$. 
Fix a rational prime $p \nmid 2 N$ and an embedding $\overline{\bf{Q}} \rightarrow \overline{\bf{Q}}_p$. Assume  
that $f$ is $p$-ordinary in the sense that the image of its $T_p$-eigenvalue under this embedding is a $p$-adic unit. 
Let $K$ be an imaginary quadratic field of discriminant $D < 0$ prime to $p$. 
Fix a $\mathcal{W}_0 = \rho_0 $ a ring class character factoring through the torsion subgroup $G_0 = \operatorname{Gal}(R_{\infty}/K)_{\operatorname{tors}}$. 
Let $k =0$ or $1$ according as to whether the pair $(f, \mathcal{W})$ is generic or exceptional respectively.  

\begin{itemize}

\item[(i)] Assume that the cyclotomic $\mu$-invariant of $f$ over $K$ vanishes. 
There exists an integer $c(k) \geq k$ such that for each character $\psi$ of conductor $c \geq c(k)$ 
factoring through the cyclotomic Galois group $\Gamma = \operatorname{Gal}(K^{\operatorname{cyc}}/K) \approx {\bf{Z}}_p$, 
the central value $L^{(k)}(1/2, f \times \rho_0 \rho_w \psi)$ does not vanish for any character $\rho_w$ of $\Omega$; in other words,
$L^{(k)}(1/2, f \times \rho \psi) \neq 0$ for any ring class character $\rho = \rho_0 \rho_w$ of 
$\varOmega = \operatorname{Gal}(K[p^{\infty}]/K)$ whose restriction to $G_0$ is $\rho_0$. \\

\item[(ii)] Assume that the level $N$ is squarefree. 
Assume as well that the residual (mod $p$) Galois representation $\overline{\rho}_f$ associated to $f$ is surjective, and that for each prime divisor 
$q$ of $N$ which remains inert in $K$ and satisfies the congruence $q \equiv \pm 1 \operatorname{mod} p$, 
this residual representation $\overline{\rho}_f$ is ramified at $q$. 
There exists an integer $d(k) \geq k$ such that for each character $\rho_w$ of conductor $d \geq d(k)$ 
factoring through the anticyclotomic Galois group $\Omega = \operatorname{Gal}(D_{\infty}/K) \approx {\bf{Z}}_p$, the value 
$L^{(k)}(1/2, f \times \rho_0 \rho_w \psi)$ does not vanish for any character $\psi$ of $\Gamma$. \end{itemize} 
  
\end{theorem} 

\begin{remark}[Remarks.] (1) Here, the conditions that $f$ be of weight $2$ and level prime to $p$ can be lifted at the expense of clarity, as can the condition that $f$ 
be $p$-ordinary. In brief, similar arguments to what we give below can be invoked for the more general setting, using other constructions of the two-variable $p$-adic 
$L$-function(s) along with other self-dual nonvanishing estimates as input. (2) Although this result might seem too strong from the analytic viewpoint, it is fact predicted 
by the refined conjecture of Birch-Swinnerton-Dyer using Iwasawa module structure arguments; cf.~the relevant discussion in Coates-Fukaya-Kato-Sujatha \cite{CFKS}. 

\end{remark}  
  
Again, this result generalizes the well-known theorems of Rohrlich \cite{Ro}, \cite{Ro2}, and Vatsal \cite{Va} for the underlying one-variable settings, 
with $\mathcal{W} = \psi = \chi \circ {\bf{N}}$ a cyclotomic character in \cite{Ro}, and with $\mathcal{W} = \rho$ a ring class character in \cite{Va}. Although 
this problem has now been addressed via purely analytic methods in \cite{VO5}, the novelty here is that we establish a minimal nonvanishing criterion 
(see Proposition \ref{criterion}) from which the result is deduced from the existence of a suitable two-variable $p$-adic $L$-function granted its nontriviality at a single, 
genuine compositum character $\mathcal{W}$ of $\mathcal{G}$ (having nontrivial anticyclotomic and cyclotomic parts). It seems likely that such a reduction 
(which does not follow directly from the Weierstrass preparation theorem) applies in other arithmetic settings, and so we spell out the deduction here as a 
self-contained sequel to \cite{VO5}. 

\subsubsection*{Acknowledgements} It is a pleasure to thank Ralph Greenberg and Michael Spiess for helpful discussions, 
as well as John Coates, Christophe Cornut, Fred Diamond, Philippe Michel, Paul Nelson, Rob Pollack, David Rohrlich, 
and an anonymous referee for helpful comments and suggestions for improvement of the writing. 

\section{Nonvanishing via $p$-adic $L$-functions}

We now prove Theorem \ref{main}, starting with background on Iwasawa algebras.
 
 \subsection{Iwasawa algebras} 

Let $\mathcal{O}$ be a discrete valuation ring, and $\mathcal{G}$ a profinite group. 
Consider the $\mathcal{O}$-Iwasawa algebra $\mathcal{O}[[\mathcal{G}]]$, which is the completed group ring 
\begin{align}\label{cgr}\mathcal{O}[[\mathcal{G}]] &= \varprojlim_{\mathcal{U}} \mathcal{O}[\mathcal{G}/ \mathcal{U}]. \end{align}
Here, the projective limit ranging over all open normal subgroups $\mathcal{U}$ of $\mathcal{G}$, and each 
$\mathcal{O}[\mathcal{G}/\mathcal{U}]$ denotes the usual group ring of the finite quotient group $\mathcal{G}/\mathcal{U}$ over $\mathcal{O}$. 
If $\mathcal{G}$ is finitely generated and abelian, then the elements of $\mathcal{O}[[\mathcal{G}]]$ can be viewed in a natural way as 
$\mathcal{O}$-valued measures on $\mathcal{G}$ (see e.g.~\cite[II, $\S$ 7]{MSD}).
In short, given $\mathcal{W}$ a finite-order character of $\mathcal{G}$, and $\mathcal{L}$ an element of $\mathcal{O}[[\mathcal{G}]]$, 
we can integrate $\mathcal{W}$ against $\mathcal{L}$ as follows. Since $\mathcal{W}$ has finite order, it defines a locally constant function on 
$\mathcal{G}$, whence there exists an open subgroup $\mathcal{U} \subset \mathcal{G}$ such that $\mathcal{W}$ is constant modulo $\mathcal{U}$. Writing 
\begin{align*} \mathcal{L}_{\mathcal{U}} &= \sum_{\sigma \in \mathcal{G}/\mathcal{U}} c_{\mathcal{U}}(\sigma) \sigma \end{align*} 
to denote the image of $\mathcal{L}$ in $\mathcal{O}[\mathcal{G}/\mathcal{U}]$, with coefficients $c_{\mathcal{U}}(\sigma) \in \mathcal{O}$, 
we can then define the integral of $\mathcal{W}$ against $\mathcal{L}$ to be the finite sum 
\begin{align}\label{specialization} \int_{\mathcal{G}} \mathcal{W}(\sigma) d\mathcal{L}
(\sigma) &= \sum_{\sigma \in \mathcal{G}/\mathcal{U}} c_{\mathcal{U}}(\sigma) \mathcal{W}(\sigma). \end{align} 
Here, $d\mathcal{L}$ denotes the measure associated to $\mathcal{L}$. It is easy to check that this measure
does not depend on the choice of open subgroup $\mathcal{U} \subset \mathcal{G}$, and hence that it is 
determined uniquely by the construction. In what follows, we shall identify an element $\mathcal{L}$ of 
$\mathcal{O}[[\mathcal{G}]]$ implicitly with its measure $d\mathcal{L}$. 
We shall also write $\mathcal{W}(\mathcal{L})$ to denote the functional defined in $(\ref{specialization})$ above, 
and refer to this as the {\it{specialization of $\mathcal{L}$ to $\mathcal{W}$}}. 
It is easy to see in this description that any group-like element $g \in \mathcal{G}$ in $\mathcal{O}[[\mathcal{G}]]$ corresponds to the Dirac measure $dg$, 
that the product  $\mathcal{L}_1\mathcal{L}_2$ of two elements $\mathcal{L}_1, \mathcal{L}_2 \in \mathcal{O}[[\mathcal{G}]]$ corresponds to the convolution 
product $d\mathcal{L}_1 \star d \mathcal{L}_2$, and that the identity $\mathcal{I}$ in $\mathcal{O}[[\mathcal{G}]]$ corresponds to a constant measure. 

\subsection{Two-variable setting}

Let us now return to the setup described above with the ${\bf{Z}}_p^2$-extension of an imaginary quadratic field $K$,
with $\mathcal{G} = \operatorname{Gal}(R_{\infty}/K)$, where $R_{\infty}$ is the maximal abelian extension of $K$ unramified outside of $p$. 
Again, $G_0 = \mathcal{G}_{\operatorname{tors}}$ denotes the finite torsion subgroup of $\mathcal{G}$, 
and $G$ the Galois group $\operatorname{Gal}(K_{\infty}/K)$ of the ${\bf{Z}}_p^2$-extension $K_{\infty}/K$ 
so that $\mathcal{G} \approx G_0 \times G$ and $G \approx {\bf{Z}}_p^2$. 
Let $\mathcal{O}$ be the ring of integers of some finite extension of ${\bf{Q}}_p$ containing the number field generated 
by the coefficients of the eigenform $f$, as well as $t$-th roots of unity for $t = \vert G_0 \vert$.
Hence, we have have an identification of completed group rings $\mathcal{O}[[\mathcal{G}]] \approx \mathcal{O}[G_0][[G]]$,
along with an injection of completed groups rings 
\begin{align}\label{isomorphism} \mathcal{O}[[\mathcal{G}]] 
&\longrightarrow \bigoplus_{\mathcal{W}_0} \mathcal{O}[[G]], 
~~~ \mathcal{L} \longmapsto (\mathcal{W}_0(\mathcal{L}))_{\mathcal{W}_0}, \end{align} 
which is also a bijection if the order of $G_0$ is prime to $p$.
Here, the direct sum runs over all characters $\mathcal{W}_0$ of $G_0$, 
and $\left( \mathcal{W}_0(\mathcal{L})\right)_{\mathcal{W}_0}$ is the vector of specializations $\mathcal{W}_0(\mathcal{L})$ of 
$\mathcal{L}$ to each character $\mathcal{W}_0$. To be clear, we only specialize to the character $\mathcal{W}_0$ of $G_0$ here, 
and not to any character of $G \approx {\bf{Z}}_p^2$. Hence, the specializations $\mathcal{W}_0(\mathcal{L})$ define genuine 
elements of the completed group ring $\mathcal{O}[[G]]$ rather than values in the ring of integers $\mathcal{O}$. 
To mark this distinction, we shall write $\mathcal{L}(\mathcal{W}_0)$ rather than $\mathcal{W}_0(\mathcal{L})$ to denote these latter elements. 

\subsection{Two-variable $p$-adic $L$-functions}

Recall that we fix an embedding $\overline{{\bf{Q}}} \rightarrow \overline{\bf{Q}}_p$. 

\begin{theorem}[Hida, Perrin-Riou]\label{hpr} 

Let $M \geq 1$ be any integer. Let $f \in S_k(\Gamma_0(M))$ be any cuspidal eigenform on the congruence subgroup $\Gamma_0(M)$
of weight $k \geq 2$ which is primitive. Assume that $f$ is $p$-ordinary, i.e. that the image of its $T_p$-eigenvalue under our fixed embedding 
$\overline{\bf{Q}} \longrightarrow \overline{\bf{Q}}_p$ is a $p$-adic unit. Assume as well that $p$ does not divide $M$.
There exists an element $\mathcal{L}_p(f) = \mathcal{L}_p(f, R_{\infty}/K)$ in the $\mathcal{O}$-Iwasawa algebra 
$\mathcal{O}[[\mathcal{G}]]$ which is characterized uniquely by the following interpolation property: For each
nontrivial finite-order character $\mathcal{W}$ of $\mathcal{G}$, we have that
\begin{align}\label{interpolation} \mathcal{W} \left(\mathcal{L}_p(f) \right) 
& = \eta(f, \mathcal{W}) \cdot \frac{L(1/2, f \times \overline{\mathcal{W}})}{8 \pi \langle f, f \rangle} \in \overline{\bf{Q}}_p.\end{align} 
Here, $\eta(f, \mathcal{W})$ denotes some explicit nonvanishing term in $\overline{\bf{Q}}_p$ (as described in \cite[Theorem 2.9]{VO3} or \cite[Th\'eor\`{e}me B]{PR88}). 
In particular, the complex value $L(1/2, f \times \overline{\mathcal{W}})$ vanishes if any only if the specialization value $\mathcal{W}(\mathcal{L}_p(f))$ vanishes.\end{theorem}

\begin{proof} The result follows from the constructions of distributions of Perrin-Riou \cite[Th\'eor\`{e}me B]{PR88} and Hida \cite{Hi}, 
using the bounded linear form construction of \cite{Hi}, whose integrality is explained in \cite[Theorem 2.9]{VO3} (with $p \geq 5$ assumed for simplicity). \end{proof}  

\begin{definition} The element $\mathcal{L}_p(f) = \mathcal{L}_p(f/R_{\infty})$ in $\mathcal{O}[[\mathcal{G}]]$ is the {\it{two-variable $p$-adic $L$-function}} 
associated to $f$ in the abelian extension $R_{\infty}/K$. \end{definition}

This $p$-adic $L$-function is characterized by the specialization property $(\ref{interpolation})$ for all nontrivial finite-order characters
$\mathcal{W}$ of $\mathcal{G} = \operatorname{Gal}(R_{\infty}/K)$. To be more precise, such characters 
$\mathcal{W} = \rho \psi = \rho (\chi \circ {\bf{N}})$ have both a dihedral ``variable" $\rho$ and a cyclotomic ``variable" 
$\psi = \chi \circ {\bf{N}}$ (hence making it a ``two-variable" $p$-adic $L$-function). 
It is said to {\it{interpolate}} the values $L(1/2, f \times \overline{\mathcal{W}})$ via the formula $(\ref{interpolation})$.

\subsection{Formal power series and the Weierstrass preparation theorem}

Let $G = \operatorname{Gal}(K_{\infty}/K) \approx {\bf{Z}}_p^2$ be the Galois group of the ${\bf{Z}}_p^2$-extension of $K$. 
Let $\gamma_1$ be a topological generator of the anticyclotomic factor $\Omega \approx {\bf{Z}}_p$, 
and let $\gamma_2$ one of the cyclotomic factor $\Gamma \approx {\bf{Z}}_p$.
We then have an isomorphism \begin{align}\label{nc} 
\mathcal{O}[[G]] &\longrightarrow \mathcal{O}[[T_1, T_2]], ~~~ (\gamma_1, \gamma_2) \longmapsto (T_1 +1, T_2 +1). \end{align} 
Here, $\mathcal{O}[[T_1, T_2]]$ denotes the ring of formal power series in the indeterminates $T_1$ and $T_2$ with coefficients 
in $\mathcal{O}$. Recall that given $\mathcal{L}$ any element of the completed 
group ring $\mathcal{O}[[\mathcal{G}]]$, we write $(\mathcal{W}_0(\mathcal{L}))_{\mathcal{W}_0}$ 
to denote its image under the injective map $(\ref{isomorphism})$. We can now consider
the image of each specialization element $\mathcal{W}_0(\mathcal{L}) = \mathcal{L}(\mathcal{W}_0)$ 
under this noncanonical isomorphism $(\ref{nc})$ to obtain a two-variable
formal power series $\mathcal{L}(\mathcal{W}_0; T_1, T_2)$. Under this composition of morphisms, 
the specialization value $\mathcal{W}(\mathcal{L})$ of $\mathcal{L}$ to a finite order character 
$\mathcal{W} = \mathcal{W}_0 \mathcal{W}_w$ of $\mathcal{G} \approx G_0 \times G$ has the following precise description. 
Let us fix a character $\mathcal{W}_0$ of the finite group $G_0$. Let us then write 
$\mathcal{W}_w = \rho_w \psi_w$ to denote any finite order character of $G = \Omega \times \Gamma$, 
where $\rho_w$ factors through the anticyclotomic quotient $\Omega$, and 
$\psi_w$ factors through the cyclotomic quotient $\Gamma$. The specialization value \begin{align*} \mathcal{W} \left(\mathcal{L}\right) 
&= \mathcal{W}_w(\mathcal{L}(\mathcal{W}_0)) = \rho_w \psi_w 
\left( \mathcal{L}(\mathcal{W}_0) \right) \end{align*} then corresponds to evaluating the associated two-variable power series 
\begin{align}\label{2VPS}\mathcal{L}(\mathcal{W}_0; T_1, T_2) &= \mathcal{L}(\mathcal{W}_0; \gamma_1-1, \gamma_2-1) \end{align} 
at the uniquely-determined roots of unity $\zeta_1 = \rho_w(\gamma_1)$ and $\zeta_2 = \psi_w(\gamma_2)$. 

We now recall the Weierstrass preparation theorem, which is a classical result that applies to power series rings in one indeterminate $T$. 
Suppose that $\mathcal{O}$ is any discrete valuation ring (or more generally any complete local ring) with maximal ideal $\mathfrak{P}$. 
Let us fix a uniformizer $\varpi_{\mathfrak{P}}$ of $\mathcal{O}$. Let $\mathcal{O}[[T]]$ denote the ring of formal power 
series in the indeterminate $T$ with coefficients in $\mathcal{O}$. Recall that a polynomial $g(T)$ in the ring $\mathcal{O}[T]$ is said to be {\it{distinguished (or Weierstrass)}} 
if it takes the form \begin{align*} g(T) &= T^r + b_{r-1} T^{r-1} + \ldots + b_0 \end{align*} for some integer $r \geq 1$, 
with each coefficient $b_i$ lying in the maximal ideal $\mathfrak{P}$. 

\begin{theorem} Let $h(T) = \sum_{j \geq 0} a_j T^j$ be an element of the formal power series ring $\mathcal{O}[[T]]$. If $h(T)$ is not identically zero, then it 
can be expressed uniquely as \begin{align*} h(T) &= u(T) g(T) \varpi_{\mathfrak{P}}^{\mu}\end{align*} of some unit $u(T)$ in $\mathcal{O}[[T]]$ times some distinguished 
polynomial $g(T)$ in $\mathcal{O}[T]$ times some integer power $\mu \geq 0$ of the fixed uniformizer $\varpi_{\mathfrak{P}}$ of $\mathcal{O}$. \end{theorem}

\begin{definition} Given a nonzero element $h(T) = u(T)g(T) \varpi_{\mathfrak{P}}^{\mu}$ of a formal power series ring $\mathcal{O}[[T]]$ as above, 
we shall refer to the degree of the distinguished polynomial $g(T)$ as the {\it{Weierstrass degree $\deg_W(h)$ of $h(T)$.}} 
Note that this degree can be characterized as the least integer $j \geq 0$ for which the coefficient $a_j$ in the power series 
expansion of $h(T) = \sum_{j \geq 0} a_j T^j$ is a unit in $\mathcal{O}$. \end{definition} 

\section{Reductions and proof of Theorem \ref{main}} 

To prove the main result, we first establish the following nonvanishing criterion to reduce the problem. 
Let us write $\mathcal{L} = \mathcal{L}_p(f, R_{\infty}/K) \in \mathcal{O}[[\mathcal{G}]]$ to
denote the two-variable $p$-adic $L$-function of Theorem \ref{hpr} above. 

\begin{proposition}[Least nontriviality criterion]\label{criterion} 
Fix a character $\mathcal{W}_0$ of the torsion subgroup of $\mathcal{G} = \operatorname{Gal}(R_{\infty}/K)$, 
and let us again write $\mathcal{L}(\mathcal{W}_0; T_1, T_2)$ to denote the image of the corresponding partially specialized $p$-adic $L$-function 
$\mathcal{L}(\mathcal{W}_0) \in \mathcal{O}[[G]]$ in the two-variable power series ring $\mathcal{O}[[T_1, T_2]]$. Suppose that ($\star$) there exists a character 
$\mathcal{W} = \mathcal{W}_0 \rho_w \psi_w$ of $\mathcal{G}$ with $\rho_w \neq {\bf{1}}$ a nontrivial character of the anticyclotomic Galois group $\Omega \approx {\bf{Z}}_p$ 
and $\psi_w \neq {\bf{1}}$ a nontrivial character of the cyclotomic Galois group $\Gamma \approx {\bf{Z}}_p$ for which $L(1/2, f \times \mathcal{W}_0 \rho_w \psi_w) \neq 0$. \\

\noindent (i) If $\operatorname{ord}_{\mathfrak{P}}(\mathcal{L}(\mathcal{W}_0; 0, T_2)) = 0$, then there exists a minimal cyclotomic exponent $\beta_0$ such that for all characters 
$\psi_w$ of $\Gamma$ of order $p^{\beta}$ with $\beta \geq \beta_0$, 
$\mathcal{L}(\mathcal{W}_0, \rho_w(\gamma_1)-1, \psi_w(\gamma_2)-1) \neq 0$ for any character $\rho_w$ of $\Omega$.\\

\noindent (ii) If If $\operatorname{ord}_{\mathfrak{P}}(\mathcal{L}(\mathcal{W}_0; T_1, 0)) =0$, then there exists a minimal anticyclotomic exponent
$\alpha_0$ such that for all characters $\rho_w$ of $\Gamma$ of order $p^{\alpha}$ with $\alpha \geq \alpha_0$, 
$\mathcal{L}(\mathcal{W}_0, \rho_w(\gamma_1)-1, \psi_w(\gamma_2)-1) \neq 0$ for any character $\psi_w$ of $\Gamma$.

\end{proposition}

\begin{proof} To show (i), we expand the $p$-adic $L$-function $\mathcal{L}(\mathcal{W}_0, T_1, T_2)$ in the cyclotomic variable $T_2$:
\begin{align}\label{cycex} \mathcal{L}(\mathcal{W}_0; T_1, T_2) &= \sum_{j \geq k} a_j(T_1) T_2^j \in \mathcal{O}[[T_1]] [[T_2]]. \end{align}
Here, $k \in \lbrace 0, 1 \rbrace$ determines the anticyclotomic root number $-\omega(N) = (-1)^k$, and each coefficient
$a_j(T_1)$ is a power series in $\mathcal{O}[[T_1]]$. Taking for granted the condition ($\star$), we deduce that $\mathcal{L}(\mathcal{W}_0, T_1, T_2)$
has a (finite) well-defined Weierstrass degree in the indeterminate $T_2$. 
Assuming as we do that $\mathcal{L}(\mathcal{W}_0, 0, T_2) \not\equiv 0 \operatorname{mod} \mathfrak{P}$, 
the Weierstrass degree is equal to the minimal integer $m \geq k$ for which the coefficient $a_m(T_1)$ is a unit in the power series ring $\mathcal{O}[[T_1]]$. 
Since units never specialize to zero, we deduce the first claim (i) by inspection (considering specializations) 
of the power series expansion $(\ref{cycex})$ in the cyclotomic variable $T_2$. 

To show (ii), we expand the $p$-adic $L$-function $\mathcal{L}(\mathcal{W}_0, T_1, T_2)$ in the anticyclotomic variable $T_1$:
\begin{align}\label{acex} \mathcal{L}(\mathcal{W}_0; T_1, T_2) &= \sum_{j \geq 0} b_j(T_2) T_1^j \in \mathcal{O}[[T_2]] [[T_1]]. \end{align}
Here, each of the coefficients $b_j(T_2)$ is a power series in $\mathcal{O}[[T_2]]$. Taking for granted the condition ($\star$) again,
we deduce that $\mathcal{L}(\mathcal{W}_0, T_1, T_2)$ has a (finite) well-defined Weierstrass degree in the indeterminate $T_1$. Assuming as 
we do that $\mathcal{L}(\mathcal{W}_0, T_1, 0) \not\equiv 0 \operatorname{mod} \mathfrak{P}$, 
this Weierstrass degree is the least integer $r \geq 0$ such that $b_r(T_2)$ is 
a unit in $\mathcal{O}[[T_2]]$. We then deduce the second claim (ii) in the same way as for (i), using that units never specialize to zero to 
read it off from the power series expansion $(\ref{acex})$ in the anticyclotomic variable $T_1$. 
\end{proof}

\begin{remark} As suggested to the author by Ralph Greenberg, the same conclusion can be drawn using Monsky's topology
of zero sets for multivariable power series rings (``${\bf{Z}}_p$-flats"), \cite{Mo}. However, our approach via expansions of the two-variable
$p$-adic $L$-functions $\mathcal{L}(\mathcal{W}_0; T_1, T_2)$ streamlines this approach considerably. \end{remark}

Using the main result of the prequel work \cite[Theorem 1.3]{VO5} to supply the criterion ($\star$), 
we then prove the following unconditional refinement. 

\begin{corollary}\label{NV} Fix a ring class character $\rho_0$ of the torsion subgroup of $\mathcal{G} = \operatorname{Gal}(R_{\infty}/K)$. 
Let $f$ be a $p$-ordinary eigenform of (squarefree) level $N$ prime to $p$ and the discriminant of $K$. \\

\noindent (i) If $\operatorname{ord}_{\mathfrak{P}}(\mathcal{L}(\mathcal{W}_0; 0, T_2)) = 0$ (i.e.~if the analytic cyclotomic $\mu$-invariant vanishes), 
then there exists a minimal exponent $\beta_0$ such that for all primitive Dirichlet characters $\chi \operatorname{mod} p^{\beta}$ of order $p^{\beta}$ 
with $\beta \geq \beta_0$, $L(1/2, f \times \rho_0 \rho_w \chi \circ {\bf{N}}) \neq 0$ for any (wild) ring class character $\rho_w$ of $\Omega$.\\

\noindent (ii) If $\operatorname{ord}_{\mathfrak{P}}(\mathcal{L}(\mathcal{W}_0; T_1, 0)) =0$ (i.e.~if the analytic anticyclotomic $\mu$-invariant vanishes), 
then there exists a minimal exponent $\alpha_0$ such that for all  ring class characters characters $\rho_w$ of $\Omega$ of order $p^{\alpha}$ 
with $\alpha \geq \alpha_0$, $L(1/2, f \times \rho_0 \rho_w \chi \circ {\bf{N}}) \neq 0$ for any primitive Dirichlet character $\chi$ of $p$-power conductor. 

\end{corollary}

\begin{proof} The result of \cite[Theorem 1.3]{VO5} implies that for an arbitrary ring class character $\rho = \rho_0 \rho_w$ factoring
through $\mathcal{G} = \operatorname{Gal}(R_{\infty}/K)$, there exists a (nontrivial) primitive Dirichlet character $\chi$ of conductor $p^{\beta}$ with
$\beta \geq 4$ (large relative to the conductor of $\rho$) for which $L(1/2, f \times \rho \chi \circ {\bf{N}}) \neq 0$. This supplies the 
nonvanishing criterion ($\star$) of Proposition \ref{criterion}. The statement then follows from Proposition \ref{criterion}. \end{proof}

Finally, let us explain what is conjectured or known about the vanishing of $\mu$-invariants for conditions (i) and (ii) above. 
The condition for Corollary \ref{NV} (i) is conjectured by Greenberg to be true in some generality\footnote{If $f$ corresponds to an elliptic curve 
which does not admit a $p$-isogeny, then the cyclotomic $\mu$-invariant is expected to vanish.}, although it appears that almost no 
progress has been made on this problem to date. More is known about the anticyclotomic $\mu$-invariant for (ii) thanks to 
Pollack-Weston \cite{PW} (see also Vatsal \cite{Va2}). In particular, we prove Theorem \ref{main} by putting together Corollary \ref{NV} with the following:

\begin{theorem}[Pollack-Weston {\cite[Theorem 1.2]{PW}}] Assume that the eigenform $f$ is $p$-ordinary, and of squarefree level $N$ not dividing the discriminant of $K$.
Assume as well that the residual (mod $p$) Galois representation $\overline{\rho}_f$ associated to $f$ is surjective, and that for each prime divisor 
$q$ of $N$ which remains inert in $K$ and satisfies the congruence $q \equiv \pm 1 \operatorname{mod} p$, this residual representation $\overline{\rho}_f$ is ramified at $q$. 
Then, for each ring class character $\rho_0$ of the torsion subgroup of $\mathcal{G} = \operatorname{Gal}(R_{\infty}/K)$, the corresponding anticyclotomic $p$-adic $L$-function
$\mathcal{L}(\rho_0; T_1, 0) \in \mathcal{O}[[T_1]]$ has vanishing $\mu$-invariant: $\operatorname{ord}_{\mathfrak{P}} (\mathcal{L}(\rho_0, T_1, 0)) = 0.$ \end{theorem}

\end{document}